\documentclass[11pt]{amsart}
\usepackage{amsmath,amssymb,amsfonts,amsthm,amsopn}
\usepackage{latexsym,graphicx}
\usepackage{pb-diagram}

\newcommand{\tfa}{time-frequency analysis}

\newcommand{\stft}{short-time Fourier transform}

\newcommand{\modsp}{modulation space}

\newtheorem{theorem}{Theorem}[section]
\newtheorem{lemma}[theorem]{Lemma}
\newtheorem{corollary}[theorem]{Corollary}
\newtheorem{proposition}[theorem]{Proposition}
\newtheorem{definition}[theorem]{Definition}

\newcommand{\beqa}{\begin{eqnarray*}}
\newcommand{\eeqa}{\end{eqnarray*}}

\newcommand{\field}[1]{\mathbb{#1}}
\newcommand{\bR}{\field{R}}        
\newcommand{\bN}{\field{N}}        
\newcommand{\bZ}{\field{Z}}        
\def\G{\mathcal{G}}

\def\la{\lambda}

\def\eps{\epsilon}

\def\cS{\mathcal{S}}

\def\cG{\mathcal{G}}
\def\cM{\mathcal{M}}

\def\cC{\mathcal{C}}

\def\a{\aleph}

\def\rd{\bR^d}

\def\rdd{{\bR^{2d}}}
\def\zdd{{\bZ^{2d}}}

\def\lrd{L^2(\rd)}

\def\R{\right)}

\def\<{\left<}
\def\>{\right>}

\def\mv1{M_v^1}

\def\phas{(x,\o )}
\def\mn{(m,n)}
\def\mn'{(m',n')}

\def\o{\xi}
\def\a{\alpha}
\def\b{\beta}

\def\R{\mathbb{R}}
\def\Ren{\mathbb{R}^d}
\def\Renn{\mathbb{R}^{2d}}

\def\sch{\mathcal{S}}

\def\Fur{\mathcal{F}}

\def\Sn2{S_{2}(L^{2}(\Ren))}
\def\S1{S_{1}(L^{2}(\Ren))}
\def\sig00{\sigma_{0,0}}

\def\la{\langle}
\def\ra{\rangle}

\begin{document}
\begin{abstract}
We perform a Gabor analysis for a large class of evolution equations with constant coefficients. We show that the corresponding propagators have a very sparse Gabor matrix, displaying off-diagonal exponential decay. The results apply to hyperbolic, weakly hyperbolic and parabolic equations. Some numerical experiments are provided. 
\end{abstract}

\title{Gabor wave packets and evolution operators}

\author{Elena Cordero,  Fabio Nicola  and Luigi Rodino}
\address{Dipartimento di Matematica,
Universit\`a di Torino, via Carlo Alberto 10, 10123 Torino, Italy}
\address{Dipartimento di Scienze Matematiche,
Politecnico di Torino, corso Duca degli Abruzzi 24, 10129 Torino,
Italy}
\address{Dipartimento di Matematica,
Universit\`a di Torino, via Carlo Alberto 10, 10123 Torino, Italy}

\email{elena.cordero@unito.it}
\email{fabio.nicola@polito.it}
\email{luigi.rodino@unito.it}

\maketitle
\section{Introduction}
The harmonic analysis represents a fundamental tool for the study of partial differential equations. Beside providing explicit expression for the solutions, it appears often as a second step of the investigation, once theorems of existence and uniqueness are established by other methods. The aim is then to provide a more precise insight to the properties of the solutions, by taking care simultaneously of the values of the function in the space domain, as well as of the frequency components. This proceeding is sometimes named micro-local analysis, synonym of time-frequency analysis or phase-space analysis. \par
Ideally, one would like to know exactly the frequencies occurring at a certain point for the solution. This is however out of reach, in view of the uncertainty principle of Heisenberg. So instead, we fix a partition of the phase-space into sufficiently large subsets, split consequently the function into wave packets, and establish which wave packets are present, or dominant, in the expression. \par
Such micro-local  decomposition can be done in different ways, the choice depending on the equation and on the problem under consideration. The aim is to obtain a sparse representation of the resolvent, or propagator. Namely, fixing attention on the Cauchy problem we want that the wave packets of the initial datum are moved, at any fixed time $t\not=0$, in a well determined way, so that only a controlled number of overlappings is allowed. Sparsity is extremely important in the numerical applications, by suggesting a natural proceeding of approximation. \par
In the present paper we choose as micro-local decomposition the Gabor decomposition, corresponding geometrically to a uniform partition of the phase-space into boxes, each wave packet occupying a box, essentially. Following \cite{CNRgabor}, we shall apply the Gabor decomposition to a class of evolution equations. We shall fix here attention on parabolic equations, performing some numerical experiments. \par
We begin by recalling the definition of Gabor frame, addressing to the next Section 2 for details and notation.  

Fix a function $g\in L^2(\rd)$ and consider the time-frequency shifts
\begin{equation}\label{intro1}
\pi(\lambda)g=e^{2\pi i n x}g(x-m),\quad \lambda=(m,n)\in\Lambda,
\end{equation}
for some lattice $\Lambda\subset\rdd$. The set of function $\{\pi(\lambda)g\}_{\lambda\in\Lambda}$ is called Gabor system. If moreover there exist $A,B>0$ such that
\begin{equation}\label{intro2}
A\|f\|^2_{L^2}\leq \sum_{\lambda\in\Lambda}|\langle f,\pi(\lambda)g\rangle|^2\leq B\|f\|^2_{L^2}
\end{equation}
for every $f\in L^2(\rd)$, we say that $\{\pi(\lambda)g\}_{\lambda\in\Lambda}$ is a {\it Gabor frame}; see e.g.\ \cite{ibero13,ibero30,ibero45,grochenig}. \par
Gabor frames have found important applications in signal processing and, more generally, to several problems in Numerical Analysis, see e.g. \cite{cordero-feichtinger-luef,str06}, and the references therein. More recently, the decomposition by means Gabor frames was applied to the analysis of certain partial differential equations, in particular the  constant coefficient Schr\"odinger, wave and Klein-Gordon equations \cite{bertinoro2,bertinoro3,bertinoro12,bertinoro17,kki1,kki2,kki3,MNRTT,baoxiang0,bertinoro57,bertinoro58,bertinoro58bis}. We also refer to the survey \cite{ruz} and the monograph \cite{baoxiang}. The analysis of variable coefficients Schr\"odinger-type operators was carried out in \cite{CNG,cgnr,fio3,fio1,tataru} for smooth symbols and in \cite{CNRanalitico1,CNRanalitico2} in the analytic category; see also \cite{nicola}. \par 
The fact of the matter is that, together with the decomposition of functions, say by a Gabor frame, there is a corresponding decomposition of operators; namely a linear operator $T$  can be regarded as the infinite matrix
\begin{equation}\label{intro3}
\langle T \pi(\mu)g,\pi(\lambda)g\rangle,\quad \lambda,\mu\in\Lambda.
\end{equation}
The more this matrix is sparse, the more this representation is useful, both for theoretical and numerical purposes. \par
In the applications of evolution equations, $T$ will be the propagator of some well-posed Cauchy problem, and will belong to some class of pseudodifferential operators (PSDO), or Fourier integral operators (FIO). \par
In \cite{CNRgabor} we have shown that Gabor frames may work as appropriate tool for theoretical and numerical analysis of the Cauchy problem for a large class of partial differential equations, including hyperbolic, weakly hyperbolic and parabolic equations with constant coefficients.\par
By fixing for a moment attention on the hyperbolic case, Gabor's approach may certainly look striking, since for the corresponding solutions the analysis is limited, in the most part of the literature, to the precise location of singularities in the space variables, the treatment of the frequency components being somewhat rough. Namely, in \cite{candes,cddy,hormander} and many others, the wave packets (the H\"ormander's wave-front set) are concentrated in a neighborhood, as small as we want, of each point $x_0$ in the space variables, geometrically multiplied by a conic neighborhood of $\xi_0$ in the frequency space, providing as a whole an infinite large domain.\par
So, the Gabor's approach and H\"ormander's approach are both compatible with the uncertainty principle of Heisenberg. The information given on the solutions of the hyperbolic equations are however quite different. By the Gabor analysis, in fact, we cannot identify any more where singularities exactly are, on the other hand the information on the frequency components is much more precise. \par
As disadvantage of the Gabor analysis, we also observe that Gabor frames do not work as soon as the hyperbolic operator is allowed to have non-constant coefficients. A simple example is given by the transport equation
\[
\partial_t u-\sum_{j=1}^d a_j(x)\partial_{x_j} u=0,\quad u(0,x)=u_0(x),
\]
whose solution at a fixed time $t\not=0$ is expressed by a change of variables in $u_0(x)$. A nonlinear change of variable is well-behaved with respect to H\"ormander's wave front set \cite[Theorem 8.2.4, Vol. I]{hormander}, whereas its representation with respect to Gabor frames is not sparse, cf.\ \cite{cnr-flp,cnr-global}.\par
As advantage of the Gabor decomposition, apart from detecting the frequency components, we emphasize that the same procedure works also for weakly hyperbolic equations and parabolic equations, whose numerical analysis is usually performed in a different way. Besides, for all these equations, we have exponentially sparse representation of the propagator $T$: 
\[
|\langle T \pi(\mu)g,\pi(\lambda)g\rangle|\lesssim \exp\big(-\epsilon|\lambda-\mu|^{1/s}\big),
\]
for every $\lambda,\mu$ in the lattice $\Lambda$, and for some positive constants $s,\epsilon$. \par\medskip
The contents of the next sections is the following. In Section 2 we recall some results on Gelfand-Shilov spaces, cf.\ \cite{GS,NR}, and time-frequency representations, cf. \cite{elena07,medit,grochenig,GZ,str06,ToftGS}. Section 3 is devoted to the almost-diagonalization (sparsity) of pseudodifferential operators. Basic references here are \cite{CNRgabor,GR}, see also \cite{charly06,GL09,rochberg}. Section 4 concerns applications to evolution equations. The numerical experiments, which are new with respect to \cite{CNRgabor}, are given in 4.2, 4.3.  

 \section{Preliminaries}
\subsection{Notations}
We denote the Schwartz class by
$\sch(\Ren)$ and the space of tempered
distributions by  $\sch'(\Ren)$.  We
use the brackets  $\la f,g\ra$ to
denote the extension to $\sch '
(\Ren)\times\sch (\Ren)$ of the inner
product $\la f,g\ra=\int f(t){\overline
{g(t)}}dt$ on $L^2(\Ren)$.

We denote the Euclidean norm of $ x \in {\bR}^d $ by $ |x| = \left( x_1 ^2 + \dots +x_d
^2 \right) ^{1/2}, $ and $ \langle x \rangle = ( 1 + |x|^2 )^{1/2}.$ We set $xy=x\cdot y$ for  the scalar product on
$\Ren$, for $x,y \in\Ren$.\par
 The Fourier
transform is normalized to be ${\hat
  {f}}(\o)=\Fur f(\o)=\int
f(t)e^{-2\pi i t\o}dt$.
We define the translation and modulation operators, $T$ and $M$, by
$$
T_x f(\cdot) = f(\cdot - x) \;\;\; \mbox{ and } \;\;\;
 M_x f(\cdot) = e^{2\pi i x \cdot} f(\cdot), \;\;\; x \in {\bR}^d.
$$
For $z=(x,\xi)$ we shall also write
\[
\pi(z)f=M_{\xi} T_x f.
\]


We
shall use the notation
$A\lesssim B$ to express the inequality
$A\leq c B$ for a suitable
constant $c>0$, and  $A
\asymp B$  for the equivalence  $c^{-1}B\leq
A\leq c B$.

\subsection{Gelfand-Shilov Spaces}

Gelfand-Shilov spaces can be considered a refinement of the Schwartz class, and they turn out to be useful when a more quantitative information about regularity and decay is required. Let us recall their definition and main properties; see \cite{GS,NR} for more details and proofs.

\begin{definition} Let there be given $ s, r>0$.
The Gelfand-Shilov type space $ S^{s} _{r} (\rd) $ is defined as all functions $f\in\cS(\rd)$ 
such that 
\[
|x^\a\partial^\beta f(x)| \lesssim A^{|\a|}B^{|\beta|}(\a!)^r(\beta!)^s,\quad \a,\beta\in\bN^d.
\]
for some $A,B>0$.
\end{definition}
We observe that the space $S^{s} _{r}(\rd) $ is nontrivial if and only if $ r + s \geq 1$.  So the smallest nontrivial space with $r=s$ is provided by $S^{1/2}_{1/2}(\rd)$. Every function of the type $P(x)e^{-a|x|^2}$, with $a>0$ and $P(x)$ polynomial on $\rd$, is in the class $S^{1/2}_{1/2}(\rd)$.  We observe the trivial inclusions $S^{s_1} _{r_1}(\rd)\subset S^{s_2} _{r_2}(\rd)$ for $s_1\leq s_2$ and $r_1\leq r_2$.\par
The Fourier transform maps $S^{s} _{r}(\rd)\to  S^{r} _{s}(\rd)$. Therefore for $s=r$ the spaces $S^{s} _{s}(\rd)$ are invariant under the action of the Fourier transform.

%
\begin{theorem} \label{simetria} Assume $ s>0, r>0, s+r\geq1$. For $f\in  \cS(\rd)$, the following conditions are equivalent:
\begin{itemize}
\item[a)] $f \in  S^{s} _{r}(\rd)$ .
\item[b)]  There exist  constants $A, B>0,$ such that
$$ \| x^{\a}  f \|_{L^\infty} \lesssim A^{|\a|} (\a !) ^{r}  \quad\mbox{and}\quad  \| \o^{\b}  \hat{f} \|_{L^\infty} \lesssim  B^{|\b|} ( \b!) ^{s},\quad \a,\b\in \bN^d. $$
\item[c)]  There exist  constants $A, B>0,$ such that
$$ \| x^{\a}  f \|_{L^\infty} \lesssim A^{|\a|} (\a !) ^{r}  \quad\mbox{and}\quad  \| \partial^{\b}  f \|_{L^\infty} \lesssim  B^{|\b|} ( \b!) ^{s},\quad \a,\b\in \bN^d. $$

\item[d)] There exist  constants $h, k>0,$ such that
$$
 \|f  e^{h  |x|^{1/r}}\|_{L^\infty} < \infty \quad\mbox{and}\quad \| \hat f  e^{k |\o|^{1/s}}\|_{L^\infty} < \infty.$$
\end{itemize}
\end{theorem}
The dual spaces of $S^s_r(\rd)$  are called spaces of tempered ultra-distributions and denoted by $(S^s_r)'(\rd)$. Notice that they contain the space of tempered distribution $\cS'(\rd)$.\par
Finally a kernel theorem holds as usual (\cite{nuova,mitjagin,treves}).
\begin{theorem}\label{kernelT} There exists an isomorphism between the space of linear continuous maps $T$ from $S^s_r(\rd)$ to $(S^s_r)'(\rd)$ and $(S^s_r)'(\rdd)$, which associates to every $T$ a kernel $K_T\in (S^s_r)'(\rdd)$ such that $$\la Tu,v\ra=\la K_T, v\otimes \bar{u}\ra,\quad \forall u,v \in S^s_r(\rd).$$ $K_T$ is called the kernel of $T$. \end{theorem}

\subsection{Time-frequency representations.}  We recall the basic
definition and tools from \tfa\ and  refer the  reader to \cite{grochenig} for a complete presentation.\par
Consider a distribution $f\in\cS '(\rd)$
and a Schwartz function $g\in\cS(\rd)\setminus\{0\}$, which will be called
{\it window}.
The short-time Fourier transform (STFT) of $f$ with respect to $g$ is $V_gf (z) = \langle f, \pi (z)g\rangle
$, $z=(x,\xi)\in\rd\times\rd$. The  \stft\ is well-defined whenever  the bracket $\langle \cdot , \cdot \rangle$ makes sense for
dual pairs of function or (ultra-)distribution spaces, in particular for $f\in
\cS ' (\rd )$ and $g\in \cS (\rd )$, $f,g\in\lrd$, or $f\in
(S^s_r) ' (\rd )$ and $g\in S^s_r (\rd )$.\par
The discrete counterpart of the above time-frequency representation is given by the so-called
{\it Gabor frames}. Namely, let $\Lambda=A\zdd$  with $A\in GL(2d,\R)$ (the group of real $2d\times 2d$ invertible matrices) be a lattice
of the time-frequency plane.
As anticipated in the Introduction, the set  of
time-frequency shifts $\G(g,\Lambda)=\{\pi(\lambda)g:\
\lambda\in\Lambda\}$ for a  non-zero $g\in L^2(\rd)$ is called a
Gabor system, whereas it is called {\it Gabor frame} if \eqref{intro2} holds. In that case, then there exists a dual window $\gamma\in L^2(\rd)$, such that $\cG(\gamma,\Lambda)$ is a frame, and every $f\in L^2(\rd)$ possesses the frame expansions
 \[
 f=\sum_{\lambda\in\Lambda}\langle f,\pi(\lambda)g\rangle\pi(\lambda)\gamma=\sum_{\lambda\in\Lambda}\langle f,\pi(\lambda)\gamma\rangle \pi(\lambda)g
 \]
 with unconditional convergence in $L^2(\rd)$. \par
 We finally pass to the characterization of some function spaces in terms of STFT decay. We have first of all the following basic result (cf.\  \cite{elena07, medit,GZ,T2}): if $g\in S^s_s(\rd$), $s\geq1/2$, then
 \begin{equation}\label{zimmermann2}
   f\in S^s_s(\rd)\Longleftrightarrow |V_g(f)(z)|\lesssim \exp\big({-\epsilon |z|^{1/s}}\big)\ \mbox{for some} \,\,\epsilon>0.
 \end{equation}
 When no decay is required on $f$ we still have a characterization in the following form (\cite[Theorem 3.1]{CNRgabor}).\begin{theorem}\label{teo1}
Consider $s>0$, $r>0$, $g\in S^s_r(\rd)\setminus\{0\}$. The following properties are equivalent:\\
(i) There exists a constant $C>0$ such that
\begin{equation}\label{smoothf}
|\partial^\a f(x)|\lesssim C^{|\a|}(\a!)^s,\quad x\in\rd,\,\a\in\bN^d.
\end{equation}
\noindent
(ii)  There exists a constant $C>0$ such that
\begin{equation}\label{STFTf}
|\o^\a V_gf\phas|\lesssim C^{|\a|}(\a!)^s,\quad \phas\in\rdd,\,\a\in\bN^d.
\end{equation}
(iii)  There exists a constant $\eps>0$ such that
\begin{equation}\label{STFTeps}
|V_gf\phas|\lesssim\exp\big({-\eps|\o|^{1/s}}\big),\quad \phas\in\rdd.
\end{equation}
\end{theorem}
If the equivalent conditions \eqref{smoothf}, \eqref{STFTf}, \eqref{STFTeps} are satisfied, we will say that $f$ is a Gevrey function when $s>1$, analytic if $s=1$ and ultra-analytic when $s<1$.\par
\section{Almost diagonalization of pseudodifferential operators}\label{section4}
Now we report on some results about the almost diagonalization of pseudodifferential operators having Gevrey, analytic (\cite{GR}) and ultra-analytic (\cite{CNRgabor}) symbols $\sigma(x,\xi)$. We adopt the so-called Weyl quantization, i.e. 
\[
\sigma^w f=\sigma^w(x,D) f=\int_{\rdd} e^{2\pi i(x-y)\xi} \sigma\Big(\frac{x+y}{2},\xi\Big) f(y) \, dy\, d\xi.
\]
We want to prove off-diagonal decay estimates for the Gabor matrix $\la\sigma^w \pi(z)g,\pi(w) g\ra$, $z,w\in\rdd$. The decay rate will be related to the regularity of the symbol $\sigma$. 
The key point is the following explicit formula linking the Gabor matrix with the short-time Fourier transform of its symbol (cf. \cite[Lemma 3.1]{charly06} and \cite{CNRgabor}).\begin{lemma}\label{lemma41} Consider $s\geq1/2$, $g\in S^s_s(\rd)$. Then, for $\sigma\in (S^s_s)'(\rdd)$,
\begin{equation}\label{311}
|\la\sigma^w \pi(z)g,\pi(w) g\ra|=|V_\Phi\sigma(u,v)|,\quad z,w\in\rdd,
\end{equation}
where $u=\frac{z+w}2$ and $v=j(w-z)$, and
\begin{equation}\label{312}
|V_\Phi \sigma(u,v)|=\left|\la\sigma^w \pi\left(u-\frac12 j^{-1}(v)\right)g,\pi\left(u+\frac12 j^{-1}(v)\right) g\ra\right|,
\end{equation}
where $j(z_1,z_2)=(z_2,-z_1)$, $z_1,z_2\in\rdd$, for some $\Phi\in S^s_s(\rdd)$.\par
\end{lemma}
It follows from this result and the equivalence \eqref{smoothf} $\Leftrightarrow$ \eqref{STFTeps} above, that the following decay estimates for the Gabor matrix of $\sigma^w$ hold (\cite{CNRgabor}). Notice that we have in fact a characterization. 
\begin{theorem}\label{CGelfandPseudo}  Let $s\geq1/2$, and $g\in S^s_s(\rd)\setminus\{0\}$. Then the following properties are equivalent for $\sigma\in\cC^\infty(\rdd)$:
\par {
(i)} The symbol $\sigma$ satisfies
\begin{equation}\label{simbsmooth} |\partial^\a \sigma(z)|\lesssim C^{|\a|}(\a!)^{s}, \quad \forall\, z\in\rdd,\,\forall \a\in\bN^{2d}.\end{equation}
{(ii)} There exists $\eps>0$ such that
\begin{equation}\label{unobis2s} |\langle \sigma^w \pi(z)
g,\pi(w)g\rangle|\lesssim \exp\big({-\eps|w-z|^{1/s}}\big),\qquad \forall\,
z,w\in\rdd.
\end{equation}
\end{theorem}

A similar characterization in the descrete setting, i.e.\ for Gabor frames, is slightly subtler. 
 Indeed, we use a recent result due to Gr{\"o}chenig and  Lyubarskii in \cite{GL09}. There sufficient conditions on the lattice $\Lambda=A \bZ^2$, $A\in GL(2,\R)$, are given in order for $g=\sum_{k=0}^n c_k H_k$, with $H_k$ Hermite functions, to form
a so-called Gabor (super)frame $\G(g,\Lambda)$, i.e. a frame where a dual window $\gamma$ exists, belonging to the space $S^{1/2} _{1/2} (\R)$ (cf.\  \cite[Lemma 4.4]{GL09}).
This theory transfers to the $d$-dimensional case  by taking a tensor product $g=g_1\otimes\cdots\otimes g_d\in S^{1/2} _{1/2} (\rd)$ of windows as above, which defines a Gabor frame on the lattice $\Lambda_1\times\cdots\times\Lambda_d$ and   possesses a dual window $\gamma=\gamma_1\otimes\cdots\otimes \gamma_d$ which still belongs to $S^{1/2} _{1/2} (\rd)$. 
\begin{theorem}\label{equivdiscr-cont}  Let $\G(g,\Lambda)$ be a Gabor super-frame for $\lrd$. Consider a function  $\sigma\in\cC^\infty(\rdd)$. Then the following properties are equivalent:
\par
{(i)} There exists $\eps>0$ such that the estimate \eqref{unobis2s} holds.\par
{(ii)} There exists $\eps>0$ such that
\begin{equation}\label{unobis2discr} |\langle \sigma^w \pi(\mu)
g,\pi(\lambda)g\rangle|\lesssim \exp\big({-\eps|\lambda-\mu|^{1/s}}\big),\qquad \forall\,
\lambda,\mu\in\Lambda.
\end{equation}
\end{theorem}
The above characterizations have several applications (\cite{CNRgabor}). Here we just consider the so-called sparsity property and the continuity of pseudodifferential operators on Gelfand-Shilov spaces. \par
\begin{corollary}
Under the assumptions of Theorem \ref{equivdiscr-cont}, let the Gabor matrix
$\langle \sigma^w \pi(\mu)
g,\pi(\lambda)g\rangle$ satisfy \eqref{unobis2discr}. Then it is sparse in the following sense.
Let $a$ be any column
or row of the matrix, and let
$|a|_n$ be the $n$-largest
entry of the sequence $a$.
Then, $|a|_n$
satisfies
\[
|a|_n\leq C \displaystyle \exp\big({-\epsilon n^{1/(2ds)}}\big),\quad n\in\bN
\]
for some constants $C>0,\epsilon>0$.
\end{corollary}
The main novelty with respect to the existing literature (cf. \cite{candes, guo-labate}) is the exponential as opposed to super-polynomial decay.
\begin{corollary}\label{gsooo}
Let $s\geq 1/2$ and consider a symbol $\sigma\in \cC^\infty(\rdd)$ that satisfies
\eqref{simbsmooth}. Then the Weyl operator $\sigma^w$ is bounded on $S^s_s(\rd)$.
\end{corollary}
Similarly one obtains boundedness on modulation spaces (\cite{F1,grochenig}) with weights having exponential growth; see  \cite{CNRgabor}. 
\section{Applications to evolution equations}
Consider an operator of the form 
\begin{equation}\label{operin}
P(\partial_t,D_x)=\partial_t^m +\sum_{k=1}^{m}a_k(D_x)\partial_t^{m-k},\quad t\in\R,\ x\in\R^d,
\end{equation}
where $a_{k}(\xi)$, $1\leq k\leq m$, are polynomials.
 They may be non-homogeneous, and their degree may be arbitrary (as usual, $D_{x_j}=\frac{1}{2\pi i}\partial_{x_j}$, $j=1,\ldots,d$).\par
We deal with the forward Cauchy problem
\[
\begin{cases}
P(\partial_t,D_x)u=0,\quad (t,x)\in\R_+\times\rd\\
\partial_t^k u(0,x)=u_k(x),\quad 0\leq k\leq m-1,
\end{cases}
\]
where $u_k\in\cS(\rd)$, $0\leq k\leq m-1$.
A sufficient and necessary condition for the above Cauchy problem with Schwartz data to be well posed is given by the {\it forward Hadamard-Petrowsky condition} (\cite[Section 3.10]{rauch}):\par 
{\it There exists a constant $C>0$ such that}
\begin{equation}\label{hp}
(\tau,\xi)\in\mathbb{C}\times\rd,\quad P(i\tau,\xi)=0 \Longrightarrow {\rm Im}\, \tau\geq -C.
\end{equation}
In fact one can see (\cite[pp. 126-127]{schwartz}) that the solution is then given by
\[
u(t,x)=\sum_{k=0}^{m-1} \partial_t^k E(t,\cdot)\ast \Big(u_{m-1-k}+\sum_{j=1}^{m-k-1} a_j(D_x) u_{m-k-1-j}\Big).
\]
with $E(t,x)=\Fur^{-1}_{\xi\to x} \sigma(t,\xi)$, where $\sigma(t,\xi)$ is the unique solution to
\[
\Big(\partial^m_t+\sum_{k=1}^m a_k(\xi) \partial_t^{m-k}\Big)\sigma(t,\xi)=\delta(t)
\]
supported in $[0,+\infty)\times\rd$. The distribution $E(t,x)$ is called the {\it fundamental solution} of $P$ supported in $[0,+\infty)\times\rd$. \par
We are therefore reduced to study the corresponding Fourier multiplier
\begin{equation}\label{hpm0}
\sigma^w(t,D_x)=\sigma(t,D_x) f =\Fur^{-1} \sigma(t,\cdot)\Fur f=E(t,\cdot)\ast f.
\end{equation}
(For Fourier multipliers the Weyl and Kohn-Nirenberg quantizations give the same operator).\par 
For example, for $t\geq0$, we have $\sigma(t,\xi)=\frac{\sin(2\pi|\xi|t)}{2\pi|\xi|}$ for the
wave operator $\partial^2_t-\Delta$; $\sigma(t,\xi)=\frac{\sin(t\sqrt{4\pi^2|\xi|^2+m^2})}{\sqrt{4\pi^2|\xi|^2+m^2}}$ for the
Klein-Gordon operator $\partial^2_t-\Delta+m^2$ ($m>0$); $\sigma(t,\xi)=e^{-4\pi^2|\xi|^2 t}$ for the heat operator $\partial_t-\Delta$. In all cases, $\sigma(t,\xi)=0$ for $t<0$.\par
We want to apply Theorem \ref{CGelfandPseudo} to the symbol $\sigma(t,x,\xi)=\sigma(t,\xi)$ of the multiplier $\sigma(t,D_x)$. To this end we present a suitable refinement of the Hadamard-Petrowsky condition. \par\medskip
{\it Assume that there are constants $C>0$, $\nu\geq1$ such that}
\begin{equation}\label{hpm}
(\tau,\zeta)\in\mathbb{C}\times\mathbb{C}^d,\quad P(i\tau,\zeta)=0 \Longrightarrow {\rm Im}\, \tau\geq -C(1+|{\rm Im}\,\zeta|)^\nu.
\end{equation}
We then have the following result (\cite{CNRgabor}). 
\begin{theorem}\label{E4}
Assume $P$ satisfies \eqref{hpm} for some $C>0$, $\nu\geq1$. Then the symbol $\sigma(t,\xi)$ of the corresponding propagator $\sigma(t,D_x)$ in \eqref{hpm0} satisfies the following estimates:
\begin{equation}\label{E3}
|\partial^\alpha_{\xi} \sigma(t,\xi)|\leq C^{(t+1)|\alpha|+t} (\alpha!)^s,\quad \xi\in\rd,\ t\geq0,\quad \alpha\in\bN^d,
\end{equation}
with $s=1-1/\nu$, for a new constant $C>0$.
\end{theorem}
Observe that the hypothesis $\nu\geq 1$ in the above theorem implies $0\leq s<1$. \par
As a consequence of Theorem \ref{E4} and Theorem \ref{CGelfandPseudo} we therefore obtain our main result. 
\begin{theorem}\label{teo5.2}
Assume $P$ satisfies \eqref{hpm} for some $C>0$, $\nu\geq1$, and set  $r=\min\{2,\nu/(\nu-1)\}$. If $g\in S^{1/r}_{1/r}(\rd)$ then  $\sigma(t,D_x)$  in \eqref{hpm0} satisfies
\begin{equation}\label{hpm3} |\langle \sigma(t,D_x) \pi(z)
g,\pi(w)g\rangle|\leq C \exp\big({-\eps |w-z|^{r}}\big),\qquad \forall\,
z,w\in\rdd,
\end{equation}
for some $\epsilon>0$ and for a new constant $C>0$. 
The inequality \eqref{hpm3} holds for $t$ belonging to an arbitrary bounded subset of $[0, +\infty
)$ with $\epsilon$ and $C$ fixed. 
\end{theorem}
Again we observe that $r>1$ in \eqref{hpm3}, so that we always obtain {\it super-exponential decay}. \par
We now detail some special cases of great interest, providing some numerical experiments. 
\subsection{Hyperbolic operators}
We recall that the operator $P(\partial_t,D_x)$ is called hyperbolic with respect to $t$ if the higher order homogeneous part in the symbol does not vanish at $(1,0,\ldots,0)\in\R\times\rd$, and $P$ satisfies the forward Hadamard-Petrowsky condition \eqref{hp}. This implies that the operators $a_k(D_x)$ in \eqref{operin} must have degree $\leq k$ and $P$ has order $m$. For example, the wave and Klein-Gordon operators are hyperbolic operators. However, $P$ is {\it not} required to be strictly hyperbolic, namely the roots of the principal symbol are allowed to coincide.\par
Now, if $P(\partial_t,D_x)$ is any hyperbolic operator, we always obtain Gaussian decay in the above theorem ($r=2$ in \eqref{hpm3}), at least for windows $g\in S^{1/2}_{1/2}(\rd)$. In fact, we have the following result (\cite{CNRgabor}).  
\begin{proposition}\label{pro5.4}
Assume $P(\partial_t,D_x)$ is hyperbolic with respect to $t$. Then the condition \eqref{hpm} is satisfied with $\nu=1$ for some $C>0$, and hence
\begin{equation} |\langle \sigma(t,D_x) \pi(z)
g,\pi(w)g\rangle|\leq C \exp\big({-\eps |w-z|^{2}}\big),\qquad \forall\,
z,w\in\rdd,
\end{equation}
if $g\in S^{1/2}_{1/2}(\rd)$, for some $\epsilon>0$ and for a new constant $C>0$.
\end{proposition}
\subsection{Wave equation}
Consider the wave operator $P=\partial^2_t-\Delta$ in $\R\times\R^d$, therefore $\sigma(t,\xi)=\frac{\sin(2\pi|\xi|t)}{2\pi|\xi|}$. The above Proposition \ref{pro5.4} applies, but we can also estimate the matrix decay directly, with the involved constants made explicit, by using the explicit expression of the fundamental solution. We state the result, for simplicity, in dimension $d\leq 3$. We take  $g(x)=2^{d/4}e^{-\pi|x|^2}$ as  window function, which belongs to $S^{1/2}_{1/2}(\rd)$, and moreover $\|g\|_{L^2}=1$ (Gaussian functions minimize the Heisenberg uncertainty so that they are, generally speaking, a natural choice for wave-packet decompositions). An explicit computation (\cite{CNRgabor}) gives the estimate
\[
|\langle \sigma(t,D_x) M_{\xi} T_{x} g, M_{\xi'} T_{x'} g\rangle|\leq t e^{-\frac{\pi}{2}[|\xi'-\xi|^2+(|x'-x|-t)_+^2]},\quad x,x',\xi,\xi'\in\rd,\quad d\leq3,
\]
where $(\cdot)_+$ denote positive part. \par
\begin{figure}[h] 
    \centering
    \includegraphics[width=8cm]{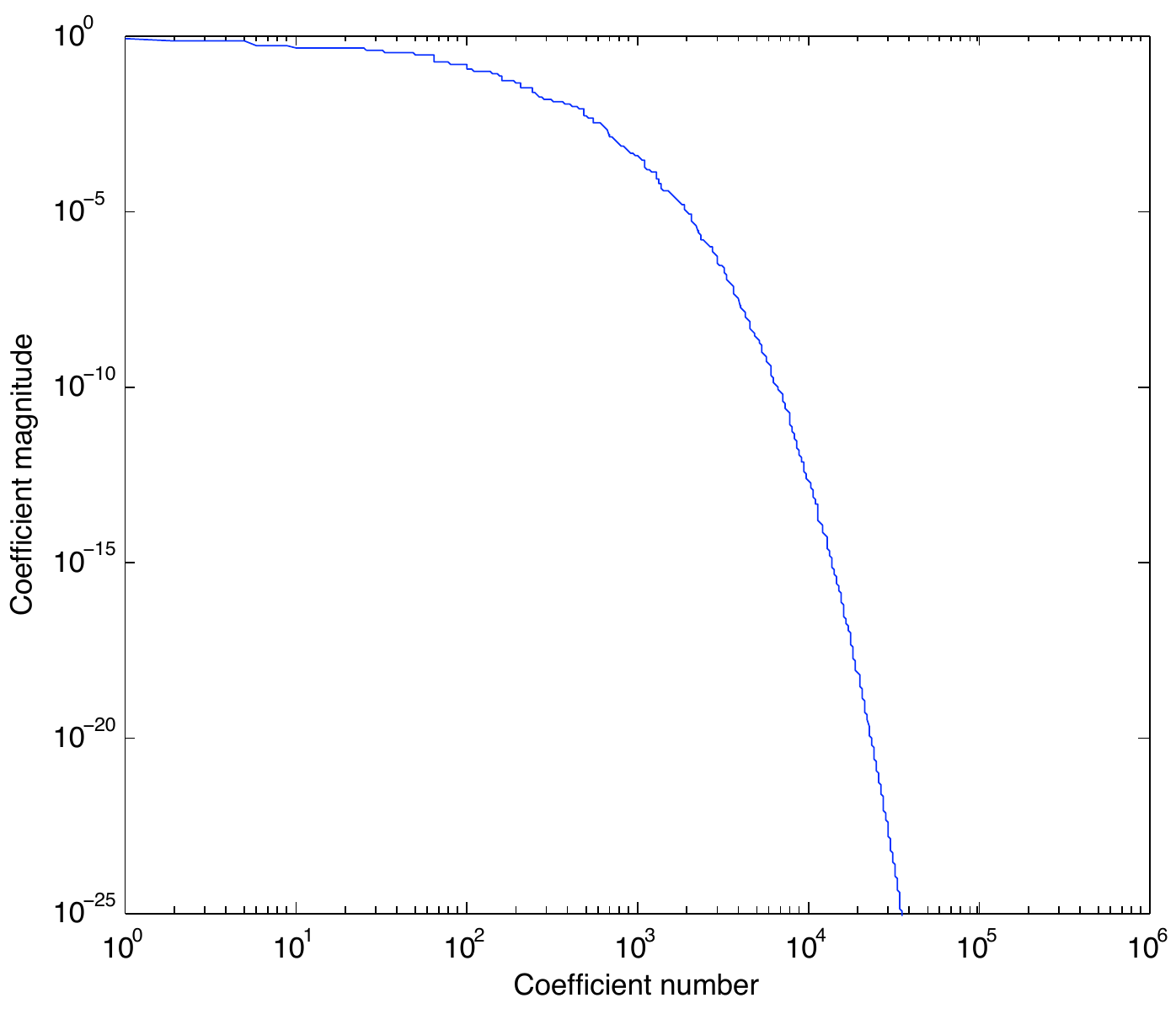}
    \caption{Decay of a generic column of the Gabor matrix for $\frac{\sin(2\pi|D| t)}{2\pi|D|}$  in dimension $d=2$ and at time $t=0.75$, with window $g(x)=\sqrt{2}e^{-\pi|x|^2}$ and lattice $\bZ^2\times(1/2)\bZ^2$.}
    \label{figura15}
 \end{figure}

Consider now the Gabor frame $\mathcal{G}(g,\Lambda)$, with $g(x)=2^{d/4}e^{-\pi|x|^2}$, $\Lambda=\bZ^d\times (1/2)\mathbb{Z}^d$ (\cite[Theorem 7.5.3]{grochenig}), and the corresponding Gabor matrix
\[
T_{m',n',m,n}=\langle \sigma(t,D_x) M_{n} T_m g, M_{n'} T_{m'} g \rangle,\quad (m,n),\ (m',n')\in\Lambda.
\]
We therefore have
\[
|T_{m',n',m,n}|\leq \tilde{T}_{m',n',m,n}:=  t e^{-\frac{\pi}{2}[|n'-n|^2+(|m'-m|-t)_+^2]},\quad (m,n),\ (m',n')\in\Lambda,\quad d\leq3.
\]
Figure \ref{figura15} shows the magnitude of the entries, rearranged in decreasing order, of a generic column, e.g. $\tilde{T}_{m',n',0,0}$ (obtained for $m=n=0$), at time $t=0.75$, in dimension $d=2$. In fact, the same figure applies to all columns, for $\tilde{T}_{m',n',m,n}=\tilde{T}_{m'-m,n'-n,0,0}$. This figure should be compared with \cite[Figure 15]{cddy}, where a similar investigation was carried out for the curvelet matrix of the wave propagator on the unit square ($d=2$) with periodic boundary conditions. It turns out that the Gabor decay is even better, in spite of the fact that we consider here the wave operator in the whole $\R^2$.

\subsection{Parabolic type equations}
Consider the operator
\begin{equation}\label{exa14}
P(\partial_t,D_x)=\partial_t+(-\Delta)^k,
\end{equation}
with $k\geq 1$ integer. In particular we get the heat operator for $k=1$. Its symbol is the polynomial
\[
P(i\tau,\zeta)=i\tau+(4\pi^2 \zeta^2)^k.
\]
An explicit computation shows that it satisfies \eqref{hpm} with $\nu=2k$. As a consequence, Theorem \ref{teo5.2} applies to $P$ with $\nu=2k$ and $r=2k/(2k-1)$.\par
In particular, the heat propagator $\sigma(t,D_x)= e^{-4\pi^2t|D|^2}$ satisfies the estimate
\[
  |\langle e^{-4\pi^2t|D|^2} \pi(z)
g,\pi(w)g\rangle|\leq C e^{-\eps |w-z|^{2}},\qquad \forall\,
z,w\in\rdd,
\]
for some $\epsilon>0$, $C>0$, if $g\in S^{1/2}_{1/2}(\rd)$. Namely, the same decay as in the case of hyperbolic equations occurs.\par
\begin{figure}[] 
    \centering
    \includegraphics[width=15cm]{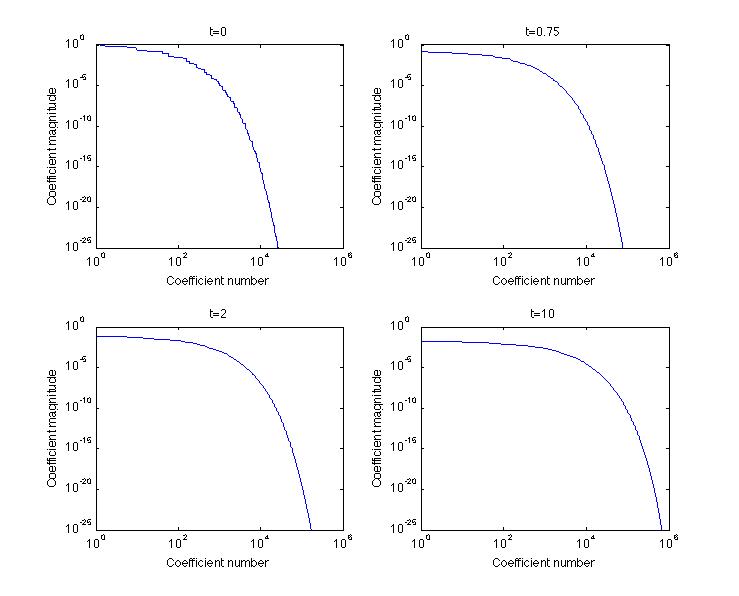}
    \caption{Decay of the column corresponding to $m=n=0$, of the Gabor matrix for the heat propagator $e^{-4\pi^2 t |D|^2}$  in dimension $d=2$ at different time instants, with window $g(x)=\sqrt{2}e^{-\pi|x|^2}$ and lattice $\bZ^2\times(1/2)\bZ^2$.}
    \label{coeff-heat}
 \end{figure}

In the following figures we summarize some numerical information about its Gabor discretization. Namely,
Figure \ref{coeff-heat} shows the decay of a column of the Gabor matrix for the heat propagator, i.e.
\[
T_{m',n',0,0}=\langle e^{-4\pi^2t|D|^2} g,M_{n'}T_{m'}g \rangle
\]
for a Gaussian window, at different time instants $t$ and in dimension $d=2$. For $t=0$ we get the identity operator, and therefore its matrix decay is the optimal one, compatibly with the uncertainty principle. As one see from the other figures the decay remains extremely good as time evolves. Also, for $t=0.75$ the decay matches that of the wave equation displayed in Figure \ref{figura15}, in spite of the fact that we no longer have here finite speed of propagation. 

\end{document}